\documentclass{article}

\usepackage{arxiv}

\usepackage[utf8]{inputenc} % allow utf-8 input
\usepackage[T1]{fontenc}    % use 8-bit T1 fonts
\usepackage{hyperref}       % hyperlinks
\usepackage{url}            % simple URL typesetting
\usepackage{booktabs}       % professional-quality tables
\usepackage{amsfonts}       % blackboard math symbols
\usepackage{nicefrac}       % compact symbols for 1/2, etc.
\usepackage{microtype}      % microtypography
\usepackage{lipsum}		% Can be removed after putting your text content
\usepackage{graphicx}
\usepackage{natbib}
\usepackage{doi}
\usepackage{amsmath}

\title{Generation and Distribution of Prime Numbers Using a Modified Lagrange Polynomial}

\author{{\hspace{1mm}Dileep Sivaraman, Branesh M. Pillai, Jackrit Suthakorn and Songpol Ongwattanakul}\thanks{Center for Biomedical and Robotics Technology (BART LAB), Faculty of Engineering, Mahidol University, Phutthamonthon District, Nakhon Pathom 73170,
Thailand, Corresponding author: Songpol Ongwattanakul (e-mail: songpol.ong@mahidol.ac.th).} 
}

\renewcommand{\shorttitle}{\textit{arXiv} Template}

\hypersetup{
pdftitle={A template for the arxiv style},
pdfsubject={q-bio.NC, q-bio.QM},
pdfauthor={David S.~Hippocampus, Elias D.~Striatum},
pdfkeywords={First keyword, Second keyword, More},
}

\begin{document}
\maketitle

\begin{abstract}
	A modified Lagrange Polynomial is introduced for polynomial extrapolation, which can be used to estimate the equally spaced values of a polynomial function. As an example of its application, this article presents a prime-generating algorithm based on a 1-degree polynomial that can generate prime numbers from consecutive primes. The algorithm is based on the condition that infinitely many prime numbers exist that satisfy the equation  $\Pi_{n} =2\Pi_{n-1} - \Pi_{n-2} \pm 2   \ \ \forall \ \Pi_{n} >7$. where $\Pi_{n-1}$ and $\Pi_{n-2}$ are the consecutive primes.
\end{abstract}

% keywords can be removed

\section{Introduction}
This study explains a new method of polynomial extrapolation called 'special addition in polynomials'(SAP), which is a modified Lagrange Polynomial that can approximate a given set of data based on a series of previous values. The salient aspect of this method is that we did not consider the values of $ x_i$ when determining the next values of f(x) $\in$ ${\rm I\!R}$. The most common method of polynomial extrapolation is Lagrange interpolation, or Newton's method of finite differences, to generate a Newton series that matches the data. Data were extrapolated using polynomials. However, these methods require values of $ x_i$ for the polynomial function f(x) \cite{waring1779vii, celant2016interpolation, wiener1964extrapolation}.

\section{The Special addition in polynomials}

Consider a polynomial function of degree m, f(x):${\rm I\!R}$ $\to$ ${\rm I\!R}$ is called a polynomial with coefficients in ${\rm I\!R}$, if there exist $a_0,a_1...a_m$ $\in$ ${\rm I\!R}$ such that \cite{barbeau2003polynomials},
\begin{equation} \label{eq1}
f(x) = a_0+a_1x+a_2x^2+...+a_mx^m\ \forall \ x\in {\rm I\!R}.
\end{equation}

\textbf{Theorem 1.} \textit{
Given a set of m+1 data points $(x_0, Y_0)$,...,$(x_i, Y_i)$,...,$(x_\rho, Y_\rho)$ where no two $x_i$ are the same, the interpolation polynomial in the Lagrange form is  \cite{waring1779vii}}

\begin{equation}
\label{e1}
L(x) = \sum_{i=1}^{\rho}Y_il_j(x)
\end{equation}

\textbf{Proposition. } \textit{
If  $Y_1, Y_2, Y_3... Y_\rho...Y_n$  are the corresponding values of equally spaced and consecutive values of the polynomial function of  $x_i$, and $Y_\rho$ be any value of f($x_\rho$) \cite{sivaraman2023nonlinear}.
}
 Then,

\begin{equation} \label{eq2}
Y_\rho = \sum_{q=1}^{p}{\frac{{(-1)^{(q+1)}(m+1)}} {q![m-(q-1)]!}Y_{(\rho-q)}} 
\end{equation}

or

\begin{equation} \label{eq3}
 Y_\rho = \frac{(m+1)}{1!}Y_{(\rho-1)} - \frac{m(m+1)}{2!}Y_{(\rho-2)} + \frac{m(m-1)(m+1)}{3!}Y_{(\rho-3)} ...
\end{equation}

Where $m$ is the degree of the polynomial. When points are consecutive and evenly spaced, the Lagrange equation \ref{e1} becomes equal to equation (\ref{eq2}).\\

\textbf{Proof:}
According equation (\ref{eq2}), for any $m$ degree polynomial,\\

\begin{equation} \label{eq4}
(x+y)^n =  \sum_{q=1}^{p}{\frac{{(-1)^{(q+1)}(m+1)}} {q![m-(q-1)]!}} (x+(y-q))^n 
\ \ \ \ \{m = n,  [m-(q-1)]>0 \}  
\end{equation}
For instance, the equation for the one-degree polynomial can be written as from (\ref{eq2})\cite{sivaraman2023nonlinear},\\

\begin{equation} \label{eq5a}
Y_\rho = 2Y_{(\rho-1)} - Y_{(\rho-2)}
\end{equation}
And a two-degree polynomial 
\begin{equation} \label{eq5}
Y_\rho = 3Y_{(\rho-1)} - 3Y_{(\rho-2)} + Y_{(\rho-3)}
\end{equation}

Using SAP, we can prove the above condition as
\begin{equation} 
\begin{split}
a(x+3)^2+b(x+3)+c &=3[a(x+2)^2+b(x+2)+c]- 3[a(x+1)^2+b(x+1)+c]+ [a(x)^2+b(x)+c] \\
\ &= a(x+3)^2+b(x+3) +c
\end{split}
\end{equation}\\
In another way,
\begin{equation} 
\begin{split}
(x+(y+3))^2 &= 3(x+(y+2))^2 - 3(x+(y+1))^2+ (x+(y))^2\\
 &= x^2+ 2x(y+3) + (y+3)^2 
\end{split}
\end{equation}
 
For a three-degree polynomial
\begin{equation} 
\begin{split}
(x+(y+4))^3 &= 4(x+(y+3))^3 - 6(x+(y+2))^3+ 4(x+(y+1))^3-(x+y)^3\\
 &=  x^3 + 3x^2(y+4) + 3x(y+4)^2 + (y+4)^3
\end{split}
\end{equation}
Similarly, we can use it for an m-degree polynomial, and prove it, as shown in \cite{sivaraman2023nonlinear}.

\section{Distribution of consecutive prime numbers}

A prime number is a positive integer and does not have any other positive integer factors than 1 and itself \cite{henderson2013dyslexia}. From this description, it is difficult to explicitly say if there are many primes or whether there are infinitely many primes.

\textbf{Definition 1. }\textit{
    A natural number $\Pi_n$ $\in$  $\mathbb{N}$ is a prime number if $\Pi_n\geq2$ and it has no positive divisors other than 1 and $\Pi_n$ \cite{weisstein2001prime, jameson2003prime}.
}

Suppose that one-degree polynomials can generate prime numbers. For prime number distribution, we can modify SAP as,

\textbf{Conjecture 1: } \textit{If $\Pi_{n-1}$, and $\Pi_{n-2}$ are consecutive prime numbers and there exist infinitely many prime numbers $\Pi_n$ such that, which satisfies the condition,}\\
\begin{equation} \label{eq8}
\Pi_{n} =2\Pi_{n-1} - \Pi_{n-2} \pm 2 \ \ \forall \ \Pi_{n} >7.
\end{equation}

\textit{i.e., there exist at least one prime number $\Pi_n$ such that, which lie in between}

$(2\Pi_{n-1} - \Pi_{n-2} + 2)$ and $(2\Pi_{n-1} - \Pi_{n-2} - 2)$

\textbf{Proof:}
    To prove that at least one prime number exists between $2\Pi_{n-1} - \Pi_{n-2} + 2$ and  $2\Pi_{n-1} - \Pi_{n-2} - 2$, where $\Pi_{n-1}$ and $\Pi_{n-2}$ are consecutive prime numbers, we can use the result known as Bertrand's postulate. 

\textbf{Definition 2. }\textit{
For any integer $n>1$, there exists at least one prime number $\Pi_{n}$ such that $n<\Pi_{n}<2n$.\cite{sondow2009ramanujan}
}

Now, let  $\Pi_{n-1}$ and $\Pi_{n-2}$ be consecutive prime numbers, such that  $\Pi_{n-2}$ $<$ $\Pi_{n-1}$  . We want to show that at least one prime number exists between $(2\Pi_{n-1} - \Pi_{n-2} + 2)$ and $(2\Pi_{n-1} - \Pi_{n-2} - 2)$.

We know that the difference between consecutive prime numbers is at least $2$. Therefore, $\Pi_{n-1}+2$, $\Pi_{n-1}+4$, ..., $\Pi_{n-2}-2$ are all composite numbers.

Consider the interval between $(2\Pi_{n-1} - \Pi_{n-2} + 2)$ and $(2\Pi_{n-1} - \Pi_{n-2} - 2)$. From Bertrand's postulate, there exists at least one prime number between $\Pi_{n-1}$ and $2\Pi_{n-1}-2$.\\
Therefore, there exists at least one prime number $\Pi_{n}$ between $(2\Pi_{n-1} - \Pi_{n-2} + 2)$ and $(2\Pi_{n-1} - \Pi_{n-2} - 2)$.\\

$2\Pi_{n-1} - \Pi_{n-2} + 2$ $<$ $\Pi_{n}$ $<$ $2\Pi_{n-1} - \Pi_{n-2} - 2$.\\

Let's consider the midpoint of this interval, which is:\\

$M = (2\Pi_{n-1} - \Pi_{n-2} + 2+2\Pi_{n-1} - \Pi_{n-2} - 2)/2=2\Pi_{n-1} - \Pi_{n-2} $\\

Notice that $M$ is an integer since $\Pi_{n-2}$ and $\Pi_{n-1}$ are prime numbers. Now, let's look at the difference between $M$ and $\Pi_{n-1}$:\\

$2\Pi_{n-1} - \Pi_{n-2}-\Pi_{n-1}=\Pi_{n-1} - \Pi_{n-2}$\\

This is exactly the difference between the two consecutive primes $\Pi_{n-2}$ and $\Pi_{n-1}$. Since this difference is at least 2, we can conclude that there are at least ($\Pi_{n-1}$- $\Pi_{n-1}$) / 2 - 1 composite numbers between $\Pi_{n-2}$ and $\Pi_{n-1}$.\\

This means that there are at least ($\Pi_{n-1}$- $\Pi_{n-1}$) / 2  + 1 numbers in the interval $(2\Pi_{n-1} - \Pi_{n-2} + 2)$ and $(2\Pi_{n-1} - \Pi_{n-2} - 2)$ that are not composite numbers, which includes prime numbers.

Therefore, there exists at least one prime number in the interval between $(2\Pi_{n-1} - \Pi_{n-2} + 2)$ and $(2\Pi_{n-1} - \Pi_{n-2} - 2)$, which proves the statement.

\textbf{Theorem 2.} \textit{
There are infinitely many primes \cite{jameson2003prime, goldston2007there}
}

\textbf{Proof:}
     Let $\Pi_{n-1}$ and $\Pi_{n-2}$ are consecutive prime numbers and according to the Conjecture-1 and equation \ref{eq8}, there exists at least one prime number $\Pi_n$ such that, which lies in between\\

($\Pi_{n} =2\Pi_{n-1} - \Pi_{n-2} + 2)\ and \ (\Pi_{n} =2\Pi_{n-1} - \Pi_{n-2} - 2$)\\

Therefore, this can be continued indefinitely and there are infinitely many primes

\section{Statistical Validation}

\begin{figure}[htbp]
    \centering
    \includegraphics[width=10cm]{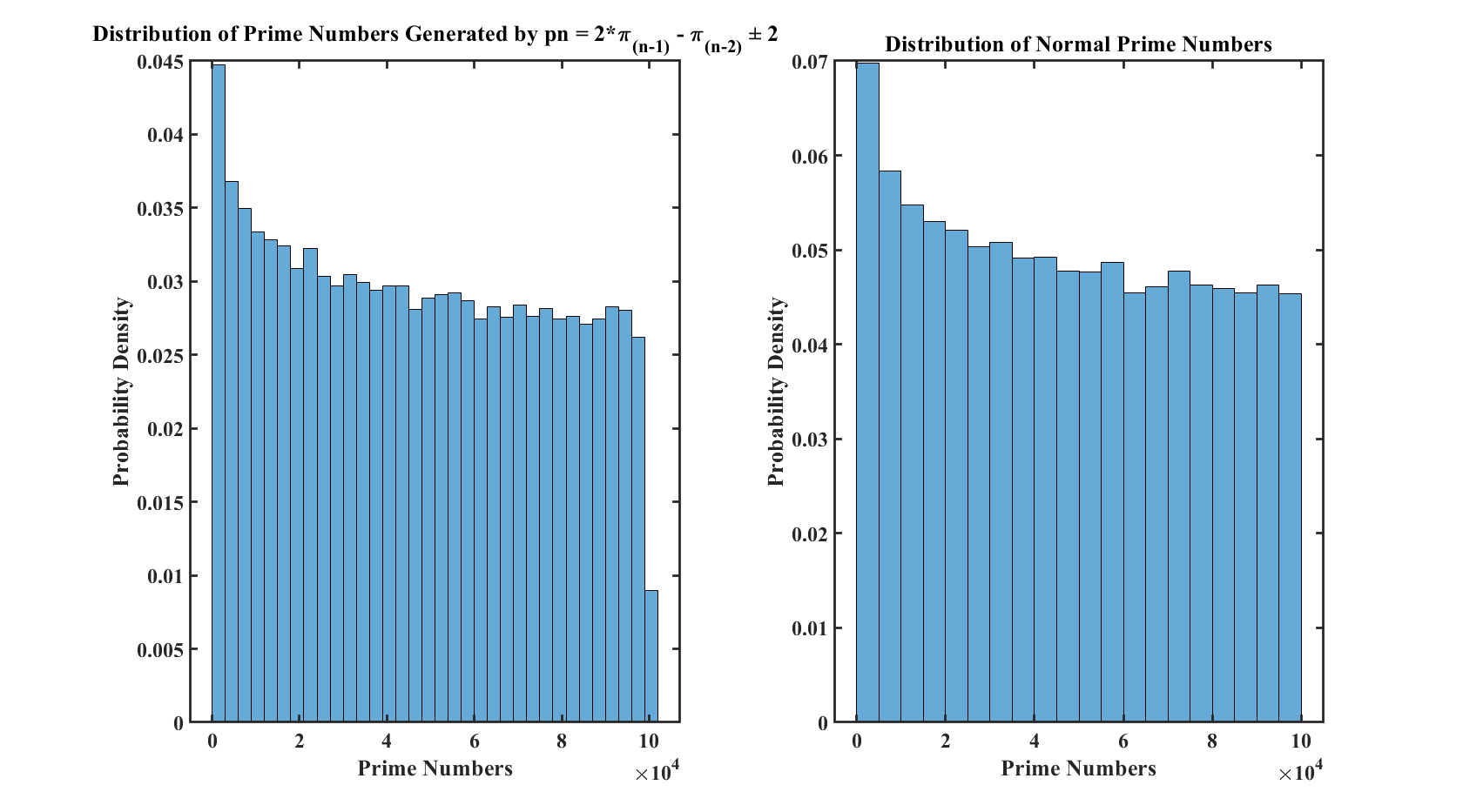}
    \caption{The histogram of the prime numbers produced by $\Pi_{n} =2\Pi_{n-1} - \Pi_{n-2} \pm 2 \ \ \forall \ \Pi_{n} >7$ appears to be different from the normal prime number distribution because the ranges of the two outputs are different. The normal prime number distribution considers all prime numbers within a certain range, whereas the prime numbers produced by $\Pi_{n}$ consider only a subset of prime numbers that satisfy a specific condition. In addition, the prime numbers produced by $\Pi_{n}$ may not be distributed uniformly, which can also affect the shape of the histogram.}
    \label{Diagram:1}
\end{figure}

\begin{figure}[htbp]
    \centering
    \includegraphics[width=10cm]{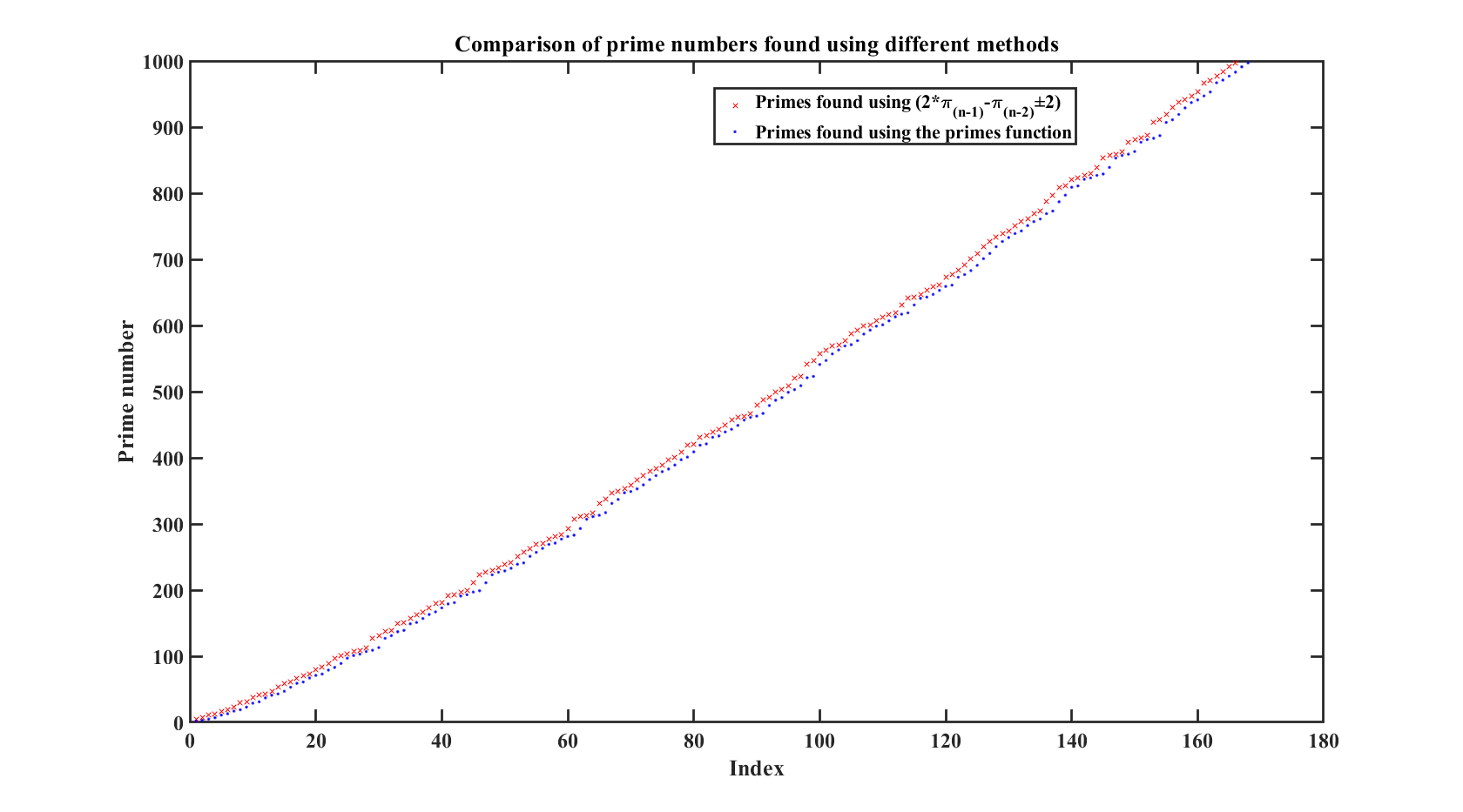}
    \caption{Using MATLAB script that finds prime numbers using a condition based on the difference between two consecutive primes. The script first initializes some variables, including the maximum number of primes to find (N = 1000). The script then checks odd numbers in the range for primality using the condition $\Pi_{n} =2\Pi_{n-1} - \Pi_{n-2} \pm 2$. The script then displays the primes found using the condition and the primes found using the primes function. The primes function is a built-in MATLAB function that finds all the prime numbers less than or equal to a given number. Finally, the script plots a graph to visualize the comparison of the two methods.}
    \label{Diagram:2}
\end{figure}

To compare the distribution of prime numbers generated by the condition  $\Pi_{n} =2\Pi_{n-1} - \Pi_{n-2} \pm 2 \ \ \forall \ \Pi_{n} >7$ with the distribution of normal prime numbers, we plotted the histograms of the two distributions and compared them visually.  MATLAB R2022b (The MathWorks, Inc., Natick, Massachusetts, United States)  generates a list of primes using the sieve of Eratosthenes \cite{weisstein2004sieve} and then loops through all pairs of consecutive primes to compute $\Pi_{n}$ and checks if they are prime. The prime numbers generated by the code are stored in the array, which is then plotted as a histogram. For comparison, we also generated a histogram of the distribution of normal prime numbers using the same number of bins. The results are presented in Figures \ref{Diagram:1} and \ref{Diagram:2}.

\section{Prime Number Spacing: Analyzing Intervals between Consecutive Primes}
It was observed that the differences between consecutive prime numbers obtained by  $\Pi_{n} =2\Pi_{n-1} - \Pi_{n-2} \pm 2$  have some pattern. This indicates that the differences between consecutive primes satisfying the condition are mostly small, whereas some outliers have relatively large gaps. The MATLAB code generates a list of primes using a sieve algorithm \cite{luo1989practical} and then checks each odd number in the range for primality. Then, a condition is applied to identify a subset of prime numbers that satisfy a specific criterion. 

\begin{figure}[htbp]
    \centering
    \includegraphics[width=10cm]{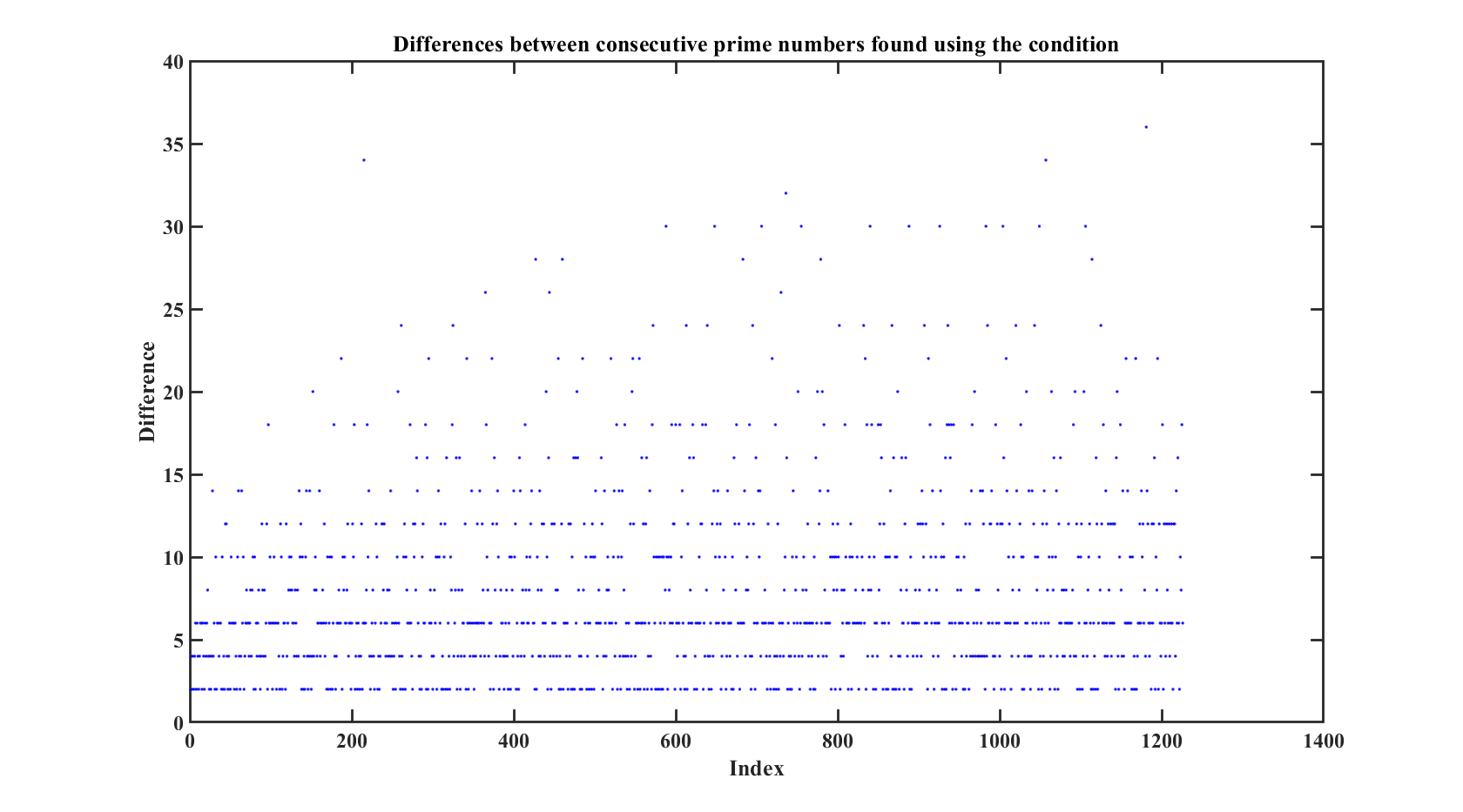}
    \caption{MATLAB utilizes a sieve algorithm \cite{luo1989practical} to generate a list of prime numbers. Subsequently, the code checks each odd number within the specified range for primality, followed by the application of a condition to identify a subset of prime numbers that satisfies a specific criterion. The distribution of the differences between consecutive primes was obtained using the condition $\Pi_{n} =2\Pi_{n-1} - \Pi_{n-2} \pm 2$, where $\Pi_{n-1}$ and $\Pi_{n-2}$ represent the two previous primes in the list. The resulting plot depicts the index of difference on the x-axis and the actual difference between consecutive primes on the y-axis. The plot revealed that the differences between consecutive primes obtained using this condition were unevenly distributed and clustered around specific values. Furthermore, the plot shows that large gaps existed between clusters of differences, signifying that the condition did not yield any primes for certain ranges. The blue dots in the plot represent the actual differences between consecutive primes.}
    \label{Diagram:3}
\end{figure}

\begin{figure}[htbp]
    \centering
    \includegraphics[width=10cm]{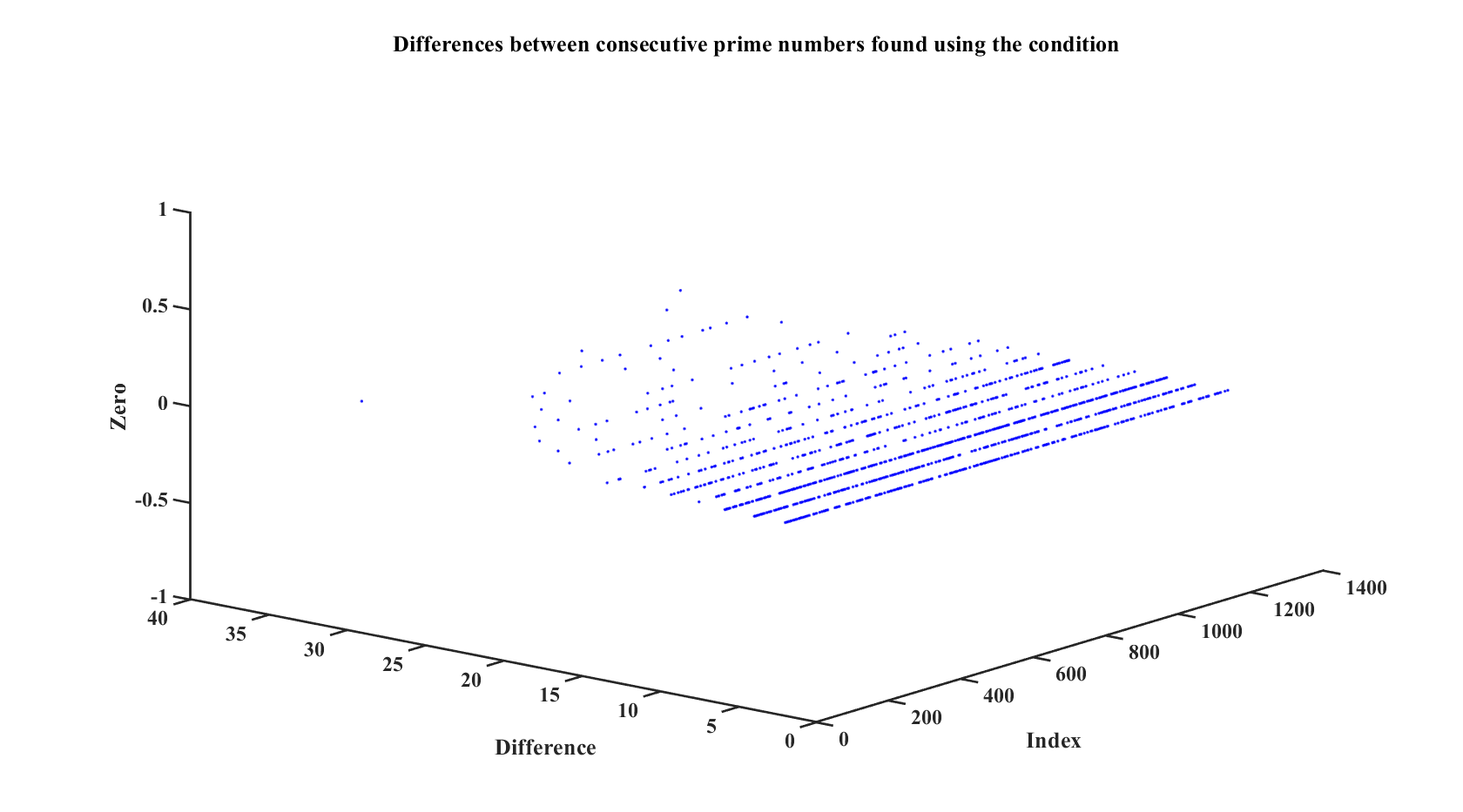}
    \caption{The graph shows the distribution of the differences between consecutive prime numbers obtained using the condition $\Pi_{n} =2\Pi_{n-1} - \Pi_{n-2} \pm 2$, where $\Pi_{n-1}$ and $\Pi_{n-2}$ are the two previous primes in the list. The x-axis represents the index of differences, the y-axis represents the actual difference between consecutive primes, and the z-axis represents zero.}
    \label{Diagram:3a}
\end{figure}

The algorithm checks odd numbers in a given range for primality using a condition and stores the primes found using the condition. It initializes variables, checks odd numbers for primality using the condition, calculates the differences between consecutive primes found using the function, and creates a 3D plot of the differences between consecutive primes. The results are presented in Figures \ref{Diagram:3} and \ref{Diagram:3a}.

\section{Prime Twins}

The prime twins were generated using the condition $\Pi_{n} =2\Pi_{n-1} - \Pi_{n-2} \pm 2$. The code loops through a range of odd numbers and checks if the two conditions are satisfied for each pair of consecutive primes. If the conditions are satisfied, then the pair of primes is added to the list of prime twins. The code loops odd numbers from 3 to 10,000, checking each number and the next largest prime number for primality using the MATLAB isprime function. If condition $\Pi_{n} =2\Pi_{n-1} - \Pi_{n-2} \pm 2$ is satisfied, the pair is added to the list, indicating that the numbers are prime twin pairs.

\begin{figure}[htbp]
    \centering
    \includegraphics[width=10cm]{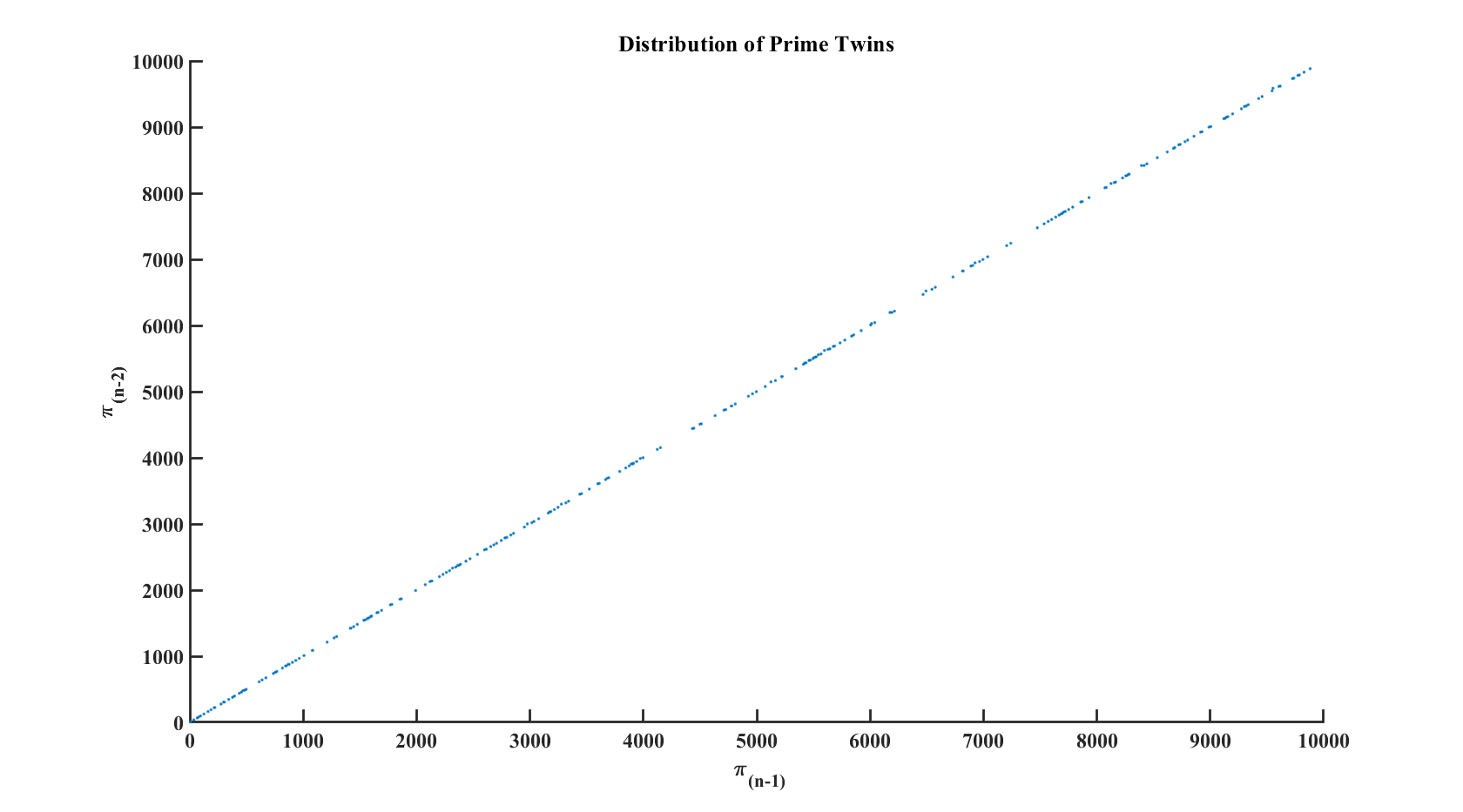}
    \caption{The code generates a list of prime twins by checking if the conditions $\Pi_{n} =2\Pi_{n-1} - \Pi_{n-2} \pm 2$ are satisfied for each pair of consecutive primes within a range of 10000 numbers. The resulting prime twin pairs are plotted in a scatter plot that shows the distribution of prime twins. The graph reveals that prime twins tend to cluster together in some areas while being more spread out in others.}
    \label{Diagram:4}
\end{figure}

The graph shows the distribution of the prime twin pairs within a range of numbers. It is interesting to observe that there are areas where prime twins seem to cluster together, whereas, in other areas, they are more widespread. The results are shown in Figure.\ref{Diagram:4}.

\begin{figure}[htbp]
    \centering
    \includegraphics[width=10cm]{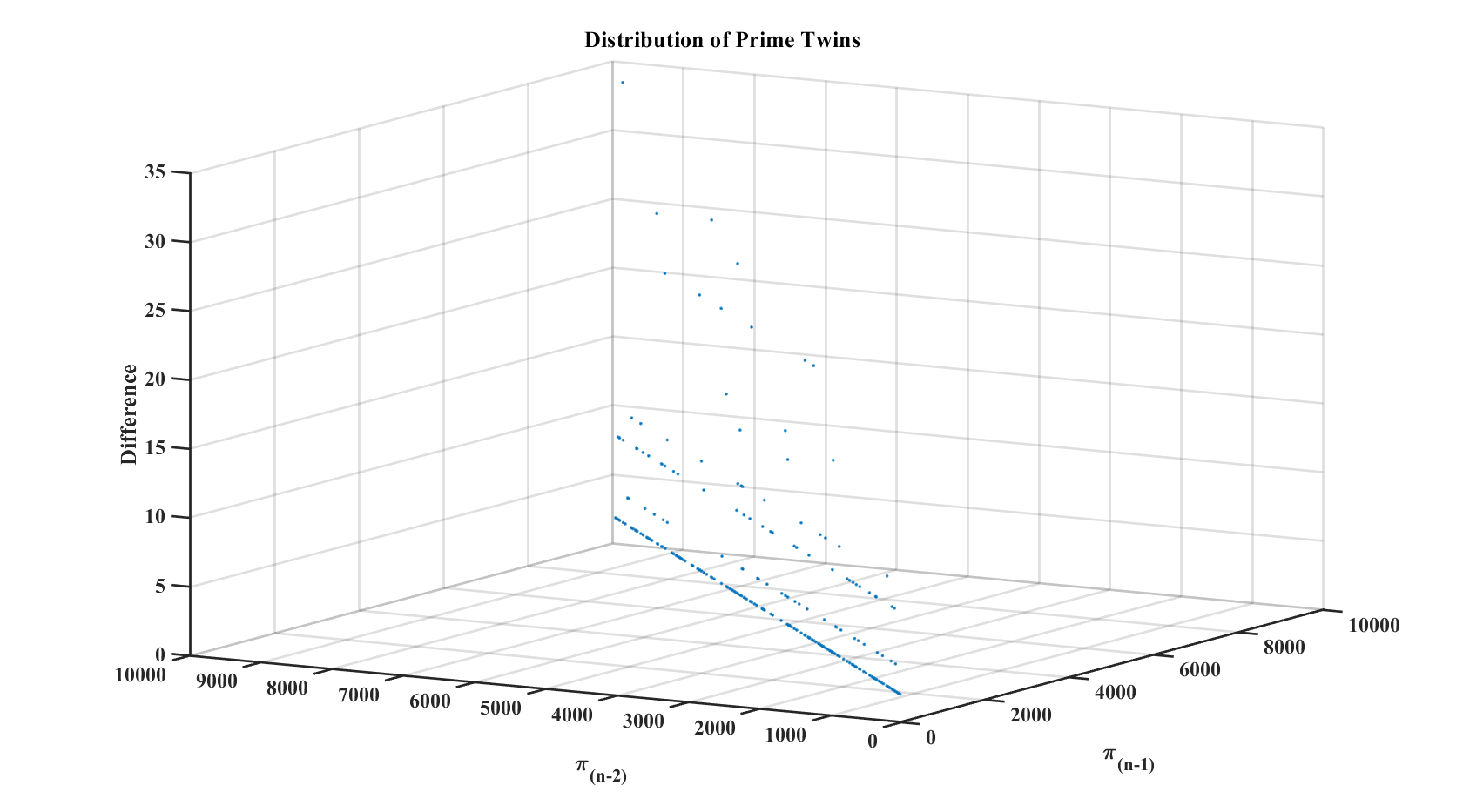}
    \caption{The x-axis represents the smaller prime number in the twin pair ($\Pi_{n-1}$), the y-axis represents the larger prime number in the twin pair ($\Pi_{n-2}$), and the z-axis represents the difference between the two primes ($\Pi_{n-1}-\Pi_{n-1}$). The scatter plot provides a visualization of the distribution of prime twin pairs in a 3D space.}
    \label{Diagram:4a}
\end{figure}

Another condition is applied to the code to generate a list of prime twin pairs using condition ($\Pi_{n} =2\Pi_{n-1} - \Pi_{n-2} \pm 2$), where $\Pi_{n-1}$ and $\Pi_{n-1}$ are consecutive prime numbers. The code checks for every prime number $\Pi_{n-1}$ in the range of 3 to 10000, and finds the next prime number $\Pi_{n-2}$. It then checks whether the conditions ($\Pi_{n} =2\Pi_{n-1} - \Pi_{n-2} \pm 2$) are prime numbers. If both conditions are satisfied, the pair ($\Pi_{n-1},\Pi_{n-2}$) is added to the list, along with their difference $\Pi_{n-1}-\Pi_{n-2}$. The results are shown in Figure.\ref{Diagram:4a}.

\bibliographystyle{unsrtnat}
\bibliography{references}  %%% Uncomment this line and comment out the ``thebibliography'' section below to use the external .bib file (using bibtex) .





% \end{thebibliography}

\end{document}